\documentclass[12pt,reqno]{amsart}
\usepackage{amsmath, amsthm, amssymb, stmaryrd}
\usepackage{hyperref}
\usepackage{enumerate}
\usepackage{url}
\usepackage{color}
\usepackage{comment}

\let\OLDthebibliography\thebibliography
\renewcommand\thebibliography[1]{
  \OLDthebibliography{#1}
  \setlength{\parskip}{0pt}
  \setlength{\itemsep}{0pt plus 0.3ex}
}

\usepackage{mathtools}

\usepackage{multirow, array}
\usepackage{placeins}

\usepackage{caption} 
\advance \topmargin by -\headheight
\advance \topmargin by -\headsep
     
\evensidemargin \oddsidemargin
\newtheorem{theorem}{\bf Theorem}[section]

\newtheorem{lemma}[theorem]{\bf Lemma}

\numberwithin{equation}{section}

\newcommand*\wrapletters[1]{\wr@pletters#1\@nil}
\def\wr@pletters#1#2\@nil{#1\allowbreak\if&#2&\else\wr@pletters#2\@nil\fi}

 \def\ge{\geqslant}

\def \deg {\mathrm{deg}}

\newcommand{\RR}{\ensuremath{\mathbb R}}

\begin{document}
\title{Expanding polynomials on sets with few products}

\author{Cosmin Pohoata}
\address{California Institute of Technology, Pasadena, CA, USA}
\email{apohoata@caltech.edu}

\thanks{}
\date{}
\begin{abstract} 
In this note, we prove that if $A$ is a finite set of real numbers such that $|AA| = K|A|$, then for every polynomial $f \in \mathbb{R}[x,y]$ we have that $|f(A,A)| = \Omega_{K,\deg f}(|A|^2)$, unless $f$ is of the form $f(x,y) = g(M(x,y))$ for some monomial $M$ and some univariate polynomial $g$. 
\end{abstract}
\maketitle
\section{Introduction}
\label{intro}
Given polynomials $f\in \RR[x]$ and $g\in \RR[x,y]$, and sets $A,B\subset \RR$, we write
\[ f(A) = \{f(a) :\ a\in A \} \ \text{ and } \ g(A,B) = \{g(a,b) :\ a\in A,\ b\in B\}. \]
That is, $g(A,B)$ is the set of distinct values that can be obtained by applying $g$ on the cartesian product $A\times B$. When $g(x,y)=x+y$ or $g(x,y)=xy$, the more convenient notation $g(A,B)=A+B$ or $g(A,B)=AB$ is generally preferred. This paper will be concerned with understanding the growth of sets such as $g(A,B)$ with respect to $|A|$ and $|B|$. We will only focus on polynomials over the reals, so our story begins with the result of Elekes and R\'onyai, who in \cite{ER00} uncovered that $|f(A,B)|$ must be asymptotically larger than $|A|$ or $|B|$, unless the polynomial $f \in \RR[x,y]$ has one of the special forms $f = h(g_1(x)+g_2(y))$ and $f = h(g_1(x)\cdot g_2(y))$, for some $h,g_1,g_2\in \RR[x]$. The current best bound for this problem is the following one by Raz, Sharir, and Solymosi \cite{RSS16}.

\bigskip

\begin{theorem} \label{th:ElekesRonyai}
Let $A, B \subset \RR$ be finite sets, and let $f\in \RR[x,y]$ be of a constant degree.
Then, unless $f = h(g_1(x)+g_2(y))$ or $f = h(g_1(x)\cdot g_2(y))$ for some $h,g_1,g_2\in \RR[x]$, we have
\[ f(A,B) = \Omega\left(\min\left\{|A|^{2/3}|B|^{2/3},|A|^2,|B|^2\right\}\right).\]
\end{theorem}

\bigskip

Theorem \ref{th:ElekesRonyai} generalizes many problems from discrete geometry and additive combinatorics, and so has many applications (for example, see \cite{RSS16} or \cite{RSZ16}). While we will not aim to give a complete background, it is important to also mention that the analogue problem has been considered over different fields instead of $\mathbb{R}$, where many interesting results are also available. See for instance \cite{Tao} for a more complete account. In particular, over finite fields it is worth pointing out the result of Vu from \cite{Vu}, who classified the two variable polynomials $f(x,y) \in \mathbb{F}_{q}[x,y]$ such that $|f(A,A)|$ is large whenever $|A+A|$ is small.

\bigskip

\begin{theorem} \label{th:Vu}
Let $A$ be a subset of $\mathbb{F}_{q}$ and let $f(x,y) \in \mathbb{F}_{q}[x,y]$ be a polynomial which cannot be written as $g(L(x,y))$ for some linear polynomial $L$ and some univariate polynomial $g$. Then,
$$\max \left\{ |A+A|, f(A,A) \right\} = \Omega\left(\min \left\{ |A|^{2/3} q^{1/3}, |A|^{3/2} q^{-1/4}\right\}\right).$$
\end{theorem}

\bigskip

Motivated by Theorem \ref{th:Vu}, Shen then considered the analogue question over the reals and proved in \cite{Shen12} the following very interesting result, which in part preceeded Theorem \ref{th:ElekesRonyai}.

\bigskip

\begin{theorem} \label{th:Shen}
Let $f\in \RR[x,y]$ be a polynomial of a constant degree that is not of the form $g(L(x,y))$ for some linear polynomial $L$ and some univariate polynomial $g$. If $A$ is a finite set of real numbers, then
\[ |A + A| |f(A, A)| = \Omega\left(|A|^{5/2}\right). \]
\end{theorem}

\bigskip

The proof of Theorem \ref{th:Shen} is in some sense a generalization of Elekes's argument from \cite{Elekes} (for the particular case when $f(x,y) = xy$), which manages to replace the spectral graph theory from the finite field case \cite{Vu} with tools from incidence geometry over the reals. In particular, Theorem \ref{th:Shen} implies that when $A \subset \mathbb{R}$ satisfies $|A+A| = O(|A|)$, we have that $|f(A,A)| = \Omega(|A|^{3/2})$ for every polynomial that is not of the form $g(L(x,y))$ for some linear polynomial $L$ and some univariate polynomial $g$. Like Theorem \ref{th:ElekesRonyai}, however, this is not optimal and it is widely believed that the exponent $3/2$ could probably be replaced with $2-\epsilon$ for every $\epsilon > 0$ in general (as it is the case for $f(x,y)=xy$; see for example \cite{Soly09} for more details).

In this paper, we address the ``dual'' problem of classifying the two variable polynomials $f(x,y) \in \mathbb{R}[x,y]$ such that $|f(A,A)|$ is large whenever $|AA|$ is small. As above, it is easy to check that there are some polynomials for which this fails. For instance, consider $A = \left\{2,2^{2},\ldots,2^{N}\right\}$ for which $|AA| = 2N-1$. If we let $f$ be a single monomial such as $f(x,y)=x^{2}y^{3}$, then it is easy to check that $f(A,A) = 5N-4$. More generally, if we choose $g$ to be a polynomial in one variable and $M(x,y)$ to be a single monomial, the $f(A,A)$ will also be usually small for $f(x,y)=g(M(x,y))$. Indeed, consider say $f(x,y) = xy+x^2y^2$; then we also have that $f(A,A) = 2N-1$. Our main result shows that $g(M(x,y))$ is the {\it{only}} real enemy, in the following strong sense.

\bigskip

\begin{theorem} \label{th:main}
Let $f\in \RR[x,y]$ be a polynomial of a constant degree that is not of the form $g(M(x,y))$ for some single monomial $M$ and some univariate polynomial $g$. If $A$ is a finite set of real numbers such that $|AA| = K|A|$, then
$$|f(A,A)| = \Omega_{K}(|A|^2).$$
\end{theorem}

\bigskip

This is also optimal up to the dependence on $K$. We prove Theorem \ref{th:main} in Section 3, after introducing the required ingredients in the upcoming Section 2. 

\bigskip

\section{Preliminaries}

\bigskip

The proof of Theorem \ref{th:main} is in some sense in spirit with the proofs of Theorem \ref{th:ElekesRonyai} and Theorem \ref{th:Shen}, but it does not rely on any incidence geometry. The main new ingredient is a quantitative version of the celebrated Schmidt subspace theorem \cite{Sch} due to Amoroso and Viada \cite{AV}.

\smallskip

\begin{theorem} \label{subspace}
Let $a_{1},\ldots,a_{n} \in K$ be nonzero elements of an algebraically closed field $K$, and let $\Gamma$ be a subgroup of $K$ of finite rank $r$. Then, the number of solutions of the equation
$$a_{1}z_{1} + \ldots + a_{n}z_{n} = 1$$
with $z_{i} \in \Gamma$ and no subsum on the left hand side vanishing is at most
$$C(n,r) := (8n)^{4n^{4}(n+nr+1)}.$$
\end{theorem}

\smallskip

Schmidt's subspace theorem (together with its different variants) represents a powerful result in number theory, particularly famous for its many applications in diophantine approximation and complexity of algebraic numbers. We will not remind them here, since many excellent surveys have been written about it, so we refer the reader for instance to \cite{Bilu} and \cite{SS}. In fact, Theorem \ref{subspace} has already manifested itself in additive combinatorics as well in \cite{Chang}, where Chang noticed that one can use it to prove that $|AA| = O(|A|)$ implies $|A+A| = \Omega(|A|^{2})$. Theorem \ref{th:main} can therefore also be seen as a generalization of this phenomenon.

The next ingredient is a multiplicative version of a somewhat more unusual version of Freiman's theorem from additive combinatorics, which is essentially a combination of \cite{Chang2} and Freiman's Lemma \cite{Freiman}. See \cite{GT06} for more details. 

\smallskip

\begin{theorem} \label{FreimanM}
Let $t \geq 1$ be an integer, let $\epsilon > 0$, and let $A$ be a finite set of real numbers with $|AA| = K|A|$ and $|A| \geq CK^{2}/\epsilon$ for some absolute constant $C > 0$. Then, $A$ is a subset of a set $G$, which is of the form\footnote{If $n$ is a positive integer, $[n]$ denotes the set $\left\{0,1,\ldots,n-1\right\}$.}
$$G:= g_{1}^{[H_{1}]} \cdot \ldots \cdot g_{r}^{[H_{r}]} =\left\{ \prod_{i=1}^{r}g_{i}^{\mu_{i}}:\ \mu_{i} \in \mathbb{Z}, \mu_{i} \in [H_{i}]\right\},$$
where $r \leq \lfloor K - 1 + \epsilon\rfloor$, all the products in
$$G^{(t)} := \left\{\prod_{i=1}^{r} g_{i}^{\mu_i}:\ \mu_{i} \in \mathbb{Z}, \mu_{i} \in [tH_{i}]\right\}$$
are pairwise distinct, and 
$$|G| = H_{1} \cdot \ldots \cdot H_{r} \leq t^{K} \exp(CK^{2} \log^{3}{K})|B|.$$
\end{theorem}

\bigskip

The last two ingredients are more algebraic in nature. First, recall that a polynomial $f\in \RR[x,y]$ is said to be \emph{reducible} if there exist polynomials $f_1,f_2\in \RR[x,y]$ of positive degrees such that $f(x,y) =f_1(x,y) \cdot f_2(x,y)$. A polynomial that is not reducible is said to be \emph{irreducible}. Furthermore, we say that a polynomial $p\in \RR[x,y]$ is \emph{decomposable} if there exists a univariate polynomial $p_1$ of degree at least two and $p_2\in \RR[x,y]$ such that $p(x, y) = p_1(p_2(x, y))$. Similarly, a polynomial that is not decomposable is said to be \emph{indecomposable}.

We will need a consequence of a theorem of Stein \cite{Stein89}, which follows from the main result of \cite{Ayad02}. See \cite{RSS16} for more details.

\smallskip

\begin{theorem} \label{Stein}
If $f\in \RR[x,y]$ is indecomposable, then the polynomial $f(x,y) - \lambda$ is reducible for at most $\deg f$ values of $\lambda\in \RR$.
\end{theorem}

\smallskip

Last but not least, we will also need the classical B\'ezout theorem \cite{Cox}, which again we only state for real polynomials, as these are the main objects of our paper.

\smallskip

\begin{theorem} \label{Bezout}
Let $f$ and $g$ be two polynomials in $\mathbb{R}[x,y]$. If $f$ and $g$ vanish simultaneously on more than $(\deg f)(\deg g)$ points of $\mathbb{R}^{2}$, then $f$ and $g$ have a common non-trivial factor.
\end{theorem}

\bigskip

\section{Proof of Theorem \ref{th:main}}

\bigskip

Let $f\in \RR[x,y]$ be a polynomial that is not of the form $g(M(x,y))$ for some single monomial $M$ and some univariate polynomial $g$, and let $d$ be the degree of $f$. We will prove that
$$|f(A,A)| = \Omega_{d,K}(|A|^{2})$$
whenever $A \subset \mathbb{R}$ satisfies $|AA|=K|A|$. The dependence on $d$ and $K$ is going to be explicit, but since it is not a priority from time to time we will reserve the right to hide certain expressions under the asymptotic notation whenever it is more convenient.

First, recall that if $f$ is decomposable, then there exist a univariate $f_1$ of degree at least two and $f_2\in \RR[x,y]$ such that $f(x,y)=f_1(f_2(x,y))$. Let $(f_1,f_2)$ be a pair of such polynomials that minimizes the degree of $f_2$. In particular, this implies that $f_2$ is indecomposable. Since $f$ is of degree at most $d$, so are $f_1$ and $f_2$. Since $f_1$ is univariate, for every $a\in \RR$ there exist at most $d$ numbers $b\in \RR$ such that $f_1(b)=a$. Thus, if $|f_2(A,A)|\geq T$ holds for some positive quantity $T$, then $|f(A,A)|\ge T/d$. It then remains to derive the lower bound for the indecomposable $f_2$, which we also know it is not a single monomial $M(x,y)$ from the hypothesis. Abusing of notation, we will refer to $f_2$ as $f$ from now on, and therefore assume without loss of generality that $f$ is indecomposable and not a single monomial as well.

Next, we naturally define the following polynomial energy of $A$ by
\[ E_{f}(A) := \left|\left\{(x,y,x',y') \in A^4:\ f(x,y) = f(x',y')\right\}\right|.\]
For each $\alpha \in f(A,A)$, we also let $m_{A}(\alpha)$ denote the number of pairs $(x,y) \in A \times A$ such that $f(x,y) = \alpha$. In particular,
$$m_{A}(\alpha)^{2} = \left|\left\{(x,y,x',y') \in A^4:\ f(x,y) = f(x',y') = \alpha \right\}\right|,$$
so by Cauchy-Schwarz,
$$E_{f}(A) = \sum_{\alpha \in f(A,A)} m_{A}(\alpha)^{2} \geq \frac{1}{|f(A,A)|} \cdot \left(\sum_{\alpha \in f(A,A)} m_{A}(\alpha)\right)^{2} = \frac{|A|^{4}}{|f(A,A)|}.$$
In order to prove that $|f(A,A)| = \Omega_{d,K}(|A|^{2})$, it thus suffices to show that $E_{f}(A) = O_{d,K}(|A|^{2})$ instead. To achieve this, we will show that for most values of $\alpha \in f(A,A)$, the number of solutions in $A \times A$ to the equation $f(x,y) = \alpha$ is at most a constant which depends solely on $d$ and $K$. More precisely, we claim that
\begin{equation} \label{key}
\left| \left\{\alpha \in f(A,A):\ m_{A}(\alpha) > C\left({d +2 \choose 2},K\right) + d^{2}2^{{d+2\choose 2}} \right\} \right| = O_{d}(1),
\end{equation}
where $C\left({d+2 \choose 2},K\right)$ is the explicit constant from Theorem \ref{subspace}.

Let us first check that this claim implies that $E_{f}(A) = O_{d,K}(|A|^{2})$. For convenience, let 
$$\Upsilon(A) :=  \left\{\alpha \in f(A,A):\ m_{A}(\alpha) > C\left({d +2 \choose 2},K\right) + d^{2}2^{{d+2\choose 2}} \right\},$$
and write
\begin{equation} \label{energy}
E_{f}(A) = \sum_{\alpha \in \Upsilon(A)} m_{A}(\alpha)^{2} + \sum_{\alpha \in f(A,A) \backslash \Upsilon(A)} m_{A}(\alpha)^{2}.
\end{equation}
For every $\alpha \in f(A,A)$, it is easy to see that $m_{A}(\alpha)^{2} = O_{d}(|A|^{2})$. Indeed, recall that this quantity is the number of solutions in $A^{4}$ to $f(x,y) = f(x',y') = \alpha$, so once $x$ and $x'$ are chosen in $A$, there are at most $d$ value for each $y$ and $y'$ that can satisfy the equality. In particular, if $|\Upsilon(A)| = O_{d}(1)$, this implies that the first term in \eqref{energy} satisfies
$$\sum_{\alpha \in \Upsilon(A)} m_{A}(\alpha)^{2} = O_{d}(|A|^{2}).$$
For the second term, note that if $\alpha \not \in \Upsilon(A)$ then 
$$M:= \max_{\alpha \in f(A,A) \backslash \Upsilon(A)} m_{A}(\alpha) \leq C\left({d +2 \choose 2},K\right) + d^{2}2^{{d+2\choose 2}}= O_{d,K}(1),$$ 
therefore
$$\sum_{\alpha \in f(A,A) \backslash \Upsilon(A)} m_{A}(\alpha)^{2} \leq M \sum_{\alpha \in f(A,A) \backslash \Upsilon(A)} m_{A}(\alpha) \leq M|A|^{2} = O_{d,K}(|A|^{2}).$$
Putting these two estimates together, we indeed get that $E_{f}(A) = O_{d,K}(|A|^{2})$. We are now left to prove \eqref{key}, which will require the tools from Section 2. 

Recall that $A$ satisfies $|AA| = K|A|$. If the size of $A$ is upper bounded by a constant in terms of $K$, then there is nothing to prove since $|f(A,A)| = \Omega_{d,K}(|A|^{2})$ is trivially true, so we can safely apply Theorem \ref{FreimanM} with $\epsilon = 1$ and $t=d$. This implies that $A$ is a subset of a set $G$, which is of the form
$$G:= g_{1}^{[H_{1}]} \cdot \ldots \cdot g_{r}^{[H_{r}]} =\left\{ \prod_{i=1}^{r}g_{i}^{\mu_{i}}:\ \mu_{i} \in \mathbb{Z}, \mu_{i} \in [H_{i}]\right\},$$
where $r \leq \lfloor K \rfloor \leq K$ and all the products all the products in
$$G^{(d)} := \left\{\prod_{i=1}^{r} g_{i}^{\mu_i}:\ \mu_{i} \in \mathbb{Z}, \mu_{i} \in [dH_{i}]\right\}$$
are pairwise distinct. We also have a quantitative estimate for $|G|$, but it is not required. 

For each $\alpha \in f(A,A)$, we now analyze the number of solutions in $A \times A$ to $f(x,y)=\alpha$. Write $f$ explicitly as
$$f(x,y) := \sum_{(i,j) \in S} a_{i,j} x^{i} y^{j},$$
where $S$ is some subset of the set of pairs $\left\{(i,j):\ i,j \geq 0, i+j \leq d\right\}$ and $a_{i,j}$ is a real coefficient for each $(i,j) \in S$. 

We begin with a first key lemma.

\bigskip

\begin{lemma} \label{lem}
For every $\alpha \in f(A,A)$, the number of solutions in $A \times A$ to
$$\sum_{(i,j) \in S} a_{i,j} x^{i} y^{j} = \alpha$$
with no subsum on the left hand side vanishing is at most $C\left({d + 2 \choose 2},K\right)$. 
\end{lemma}

\smallskip

\begin{proof} Let $\Gamma$ be the multiplicative subgroup of $\mathbb{C}^{*}$ generated by $g_{1},\ldots,g_{r}$, which has rank $r \leq K$ and contains $G$ (and thus also $A$). The number of solutions to
\begin{equation} \label{sub}
\sum_{(i,j) \in S} a_{i,j} z_{i,j} = \alpha
\end{equation}
with $z_{i,j} \in G^{(d)}$ for each $(i,j) \in S$ and no subsum on the left hand side vanishing is at most the number of solutions to \eqref{sub} with the $z_{i,j}$ in $\Gamma$ and no subsum on the left hand side vanishing, so by Theorem \ref{subspace} it is at most $C\left({d + 2 \choose 2},K\right)$. If we also can argue that for each such solution $(z_{i,j})_{(i,j) \in S}$ to \eqref{sub}, there is at most one solution $(x,y) \in G \times G$ (and thus in $A \times A$) to the system of equations
\begin{equation} \label{system}
x^{i} y^{j} = z_{i,j}\ \ \ \ \text{for each}\ (i,j) \in S,
\end{equation}
then the claim follows.

For $x$ and $y$ in $G$, write
$$x = g_{1}^{x_{1}} \cdot \ldots \cdot g_{r}^{x_{r}}\ \ \text{and}\ \ y = g_{1}^{y_{1}} \cdot \ldots \cdot g_{r}^{y_{r}},$$
where $x_{k},y_{k} \in [H_{k}]$ for each $k \in \left\{1,\ldots,r\right\}$. Similarly, for $z_{i,j} \in G^{(d)}$, let
$$z_{i,j} =g_{1}^{z_{i,j,1}} \cdot \ldots \cdot g_{r}^{z_{i,j,r}},$$
where $z_{i,j,k} \in [dH_{k}]$ for each $k \in \left\{1,\ldots,r\right\}$. Plugging these expressions into \eqref{system}, we get
$$g_{1}^{ix_{1}+jy_{1}} \cdot \ldots \cdot g_{r}^{ix_{r}+jy_{r}} = g_{1}^{z_{i,j,1}} \cdot \ldots \cdot g_{r}^{z_{i,j,r}} \ \ \ \ \text{for each}\ (i,j) \in S.$$
Furthermore, since $i+j \leq d$, we also have that $ix_{k}+jy_{k} \in [dH_{k}]$ for each $k$, so by the fact that $G^{(d)}$ has all its products pairwise distinct, it follows that \eqref{system} translates into the following system of equalities, call it $\mathcal{S}_{i,j}$, satisfied by the exponents above for each $(i,j) \in S$:
$$ix_{k} + jy_{k} = z_{i,j,k}\ \ \ \text{for all}\ k \in \left\{1,\ldots,r\right\}.$$
At this point, recall that $f$ is indecomposable by our assumption and is also not a single monomial, so it must contain at least two monomials, say $x^{i}y^{j}$ and $x^{i'}y^{j'}$, for which the two-dimensional vectors $(i,j)$ and $(i',j')$ are not a (rational) scalar multiple of each other. In particular, if a pair $(x,y) \in A \times A \subset G \times G$ exists to satisfy both $\mathcal{S}_{i,j}$ and $\mathcal{S}_{i',j'}$, then each pair $(x_{k},y_{k})$ is uniquely determined in terms of $i,j,i',j'$ and $z_{i,j,k}$, $z_{i',j',k}$ for each $k \in \left\{1,\ldots,r\right\}$, which implies that $(x,y)$ is then uniquely determined. This proves the claim. 
\end{proof}

\smallskip

We now analyze what happens if there are vanishing subsums on the left hand side of $f(x,y) = \alpha$. In this sense, we prove the following second key lemma.

\bigskip

\begin{lemma} \label{steinbez}
For all but possibly at most $d+1$ values of $\alpha \in f(A,A)$, the number of pairs $(x,y) \in A \times A$ satisfying $f(x,y) = \alpha$ with some vanishing subsum on the left hand side is at most $d^{2}2^{{d+2\choose 2}}$.
\end{lemma}

\smallskip

\begin{proof}
Recall $f(x,y) := \sum_{(i,j) \in S} a_{i,j} x^{i} y^{j}$ with $|S| \geq 2$, and now suppose that 
$$\sum_{(i’,j’) \in S’} a_{i’,j’} x^{i’} y^{j’} = 0$$
for some nontrivial subset $S' \subset S$. Let $N_{S'}(\alpha)$ be number of common solutions in $A \times A$ to 
\begin{equation} \label{bez}
f(x,y)-\alpha = 0\ \ \text{and}\ \ g_{S'}(x,y) = 0,
\end{equation}
where $g_{S'} \in \mathbb{R}[x,y]$ is the polynomial defined by
$$g_{S'}(x,y):= \sum_{(i’,j’) \in S'} a_{i’,j’} x^{i’} y^{j’}$$
for a nontrivial subset $S'$ of $S$. By a union bound, it suffices to prove that
$$\sum_{S' \subset S} N_{S'}(\alpha) \leq d^{2} 2^{{d+2 \choose 2}}$$
holds for all but possibly at most $d+1$ values of $\alpha \in f(A,A)$. 

For each $S' \subset S$, note that $g_{S'}$ has degree at most $d$, since $\deg f = d$. By Theorem \ref{Stein} there are at most $d$ values of $\alpha$ for which $f(x,y) - \alpha$ is reducible, and at most one value for which $f(x,y)-\alpha$ may be identical to $g_{S'}(x,y)$, for some $S' \subset S$ ($\alpha$ may be equal to the free term in $f$). For each of the other $\alpha \in f(A,A)$, we have that $N_{S'}(\alpha) \leq d^{2}$ for every proper $S' \subset S$. Indeed, if $\alpha$ is such that the polynomial $f(x,y) - \alpha$ is irreducible in $\mathbb{R}[x,y]$ and does not coincide with $g_{S'}(x,y)$, then this simply follows from Theorem \ref{Bezout}, since \eqref{bez} must have at most $d^{2}$ solutions if there is no common factor. Therefore,
$$\sum_{S' \subset S} N_{S'}(\alpha) \leq d^{2} 2^{|S|} \leq d^{2} 2^{{d + 2 \choose 2}}$$
is indeed satisfies by all $\alpha \in f(A,A)$, except for perhaps at most $d+1$ values. This completes the proof of Lemma \ref{steinbez}.
\end{proof}

\bigskip

Claim \eqref{key} now follows by combining Lemma \ref{lem} and Lemma \ref{steinbez}. Indeed, together these two imply that for all but possibly at most $d+1$ values of $\alpha \in f(A,A)$, the number of pairs $(x,y) \in A \times A$ with $f(x,y) = \alpha$ is at most $C\left({d + 2 \choose 2},K\right) + d^{2}2^{{d+2\choose 2}}$. In other words, $\left|\Upsilon(A)\right| \leq d+1$, where
$$\Upsilon(A) :=  \left\{\alpha \in f(A,A):\ m_{A}(\alpha) > C\left({d +2 \choose 2},K\right) + d^{2}2^{{d+2\choose 2}} \right\}.$$
This completes the proof of Theorem \ref{th:main}.

\bigskip

\bigskip

{\bf{Acknowledgements}}. I would like to thank Vlad Matei, Adam Sheffer and Dmitrii Zhelezov for useful discussions.


\bigskip

\begin{thebibliography}{50}

\bibitem{AV} F. Amoroso and E. Viada, Small points on subvarieties of a torus, 
\emph{Duke Mathematical Journal}, {\bf{150}}(3):407-442, 2009.

\bibitem{Ayad02}
M.\ Ayad,
Sur les polyn\^omes $f(X,Y)$ tels que $K[f]$ est int\'egralement ferm\'e dans $K[X,Y]$,
\emph{Acta Arith.} {\bf 105} (2002), 9--28.

\bibitem{Bilu} Y. Bilu. The Many Faces of the Subspace Theorem (after Adamczewski,
Bugeaud, Corvaja, Zannier...). S\'eminaire Bourbaki, Expos\'e 967, 59\`eme
ann\'ee (2006-2007), \emph{Ast\'erisque} 317 (2008), 1-38, May 2007.

\bibitem{Chang} M.-C. Chang, Sum and product of different sets, 
\emph{Contributions to Discrete
Mathematics}, {\bf{1}} (1), 2006.

\bibitem{Chang2} M.-C. Chang, Generalized arithmetical progressions and sumsets, 
\emph{Acta Math. Hungar.}, {\bf{65}} (1994), no. 4, 379-388.

\bibitem{Cox}  D. A. Cox, J. Little and D. O'Shea, \emph{Using Algebraic Geometry}, 
Springer-Verlag, 2nd Edition, Heidelberg 2005.

\bibitem{Elekes}  G. Elekes, On the number of sums and products, {\it{Acta Arith.}}, {\bf{81}} (1997) 365-367.

\bibitem{ER00} G.\ Elekes, L.\ R\'onyai, A combinatorial problem on polynomials and rational functions, \emph{J.\ Combin.\ Theory Ser.\  A} {\bf 89} (2000), 1-20.

\bibitem{Freiman} G.A. Freiman, \emph{Foundations of a Structural Theory of Set Addition},  
Kazan Gos. Ped. Inst., Kazan, 1966 (Russian). English translation in Translations of Mathematical Monographs {\bf{37}}, Amer. Math. Soc., Providence, RI, USA, 1973.

\bibitem{GT06} 
B. Green, T. Tao, 
Compressions, Convex Geometry and the Freiman-Bilu Theorem,
\emph{Q. J. Math.}, {\bf{57}} (4) (2006), 495-504.

\bibitem{RSS16}
O.\ E.\ Raz, M.\ Sharir, and J.\ Solymosi,
Polynomials vanishing on grids: The Elekes-R\'onyai problem revisited,
\emph{Amer.\ J.\ Math.} {\bf 138} (2016), 1029-1065.
%

\bibitem{RSZ16}
O. E. Raz, M. Sharir, and F. de Zeeuw, 
Polynomials vanishing on Cartesian products: The Elekes-Szab´o Theorem revisited, 
\emph{Duke Math. J.}, {\bf{165}} (2016), 3517-3566.

\bibitem{Sch} 
W. M. Schmidt, Norm form equations, \emph{Annals of Mathematics}, {\bf{96}}(3): 526-551, 1972.

\bibitem{Shen12}
C.\ Shen,
Algebraic methods in sum-product phenomena,
\emph{Israel J.\ Math.} {\bf 188} (2012), 123-130.
%

\bibitem{Soly09} 
J. Solymosi, Bounding multiplicative energy by the sumset, 
\emph{Adv. Math.}, {\bf 222} (2009), pp. 402-408.

\bibitem{SS} 
R. Schwarz, J. Solymosi, \emph{Combinatorial Applications of the Subspace Theorem}, Geometry, Structure and Randomness in Combinatorics. Scuola Normale Superiore (2014), 123-140. 

\bibitem{Stein89}
Y.\ Stein,
The total reducibility order of a polynomial in two variables,
\emph{Israel J.\ Math.} {\bf 68} (1989), 109-122.
%

\bibitem{Tao}
T. Tao, Expanding polynomials over finite fields of large characteristic, and a regularity lemma for definable sets, 
\emph{Contributions to Discrete Mathematics} {\bf{10}} (2015), 22-98.

\bibitem{Vu}
V. Vu, 
Sum-product estimates via directed expanders, 
\emph{Math. Res. Lett.}, {\bf{15}} (2008), 375-388.








\end{thebibliography}
\end{document}